\documentclass[titlepage,11pt]{article}
\usepackage{latexsym, tikz}
\usepackage{amsfonts}
\usepackage{amsmath}
\usepackage{textcomp}
\usetikzlibrary{graphs}
\usetikzlibrary{arrows}
\usetikzlibrary{decorations.markings}

\usepackage{float}

\newcommand\blackslug{\hbox{\hskip 1pt \vrule width 4pt height 8pt depth 1.5pt
        \hskip 1pt}}
\newcommand\bbox{\hfill \quad \blackslug \bigbreak}
\def\dd{\hbox{-}}
\def\cc{\hbox{-}\cdots\hbox{-}}
\def\ll{,\ldots,}

\title{Detecting a long odd hole}
\author{Maria Chudnovsky\thanks{This material is based upon work
supported in part by the U. S. Army
Research Office under grant number W911NF-16-1-0404, and by NSF grant
DMS-1763817.}\\
Princeton University, Princeton, NJ 08544
\\
\\
Alex Scott\thanks{Supported by a Leverhulme Trust Research
Fellowship.}\\
Mathematical Institute, University of Oxford, Oxford OX2 6GG, UK
\\
\\
Paul Seymour\thanks{Supported by NSF grant DMS-1800053 and AFOSR grant
A9550-19-1-0187, and partially supported by the Simons Foundation
and by the Mathematisches Forschungsinstitut Oberwolfach.}\\
Princeton University, Princeton, NJ 08544.}

\date{January 31, 2019; revised April 24, 2020}

\newtheorem{thm}{}[section]

\newcommand{\Proof}{\noindent{\bf Proof.}\ \ }

\begin{document}
\maketitle
\begin{abstract}
For each integer $\ell\ge 5$, we give a polynomial-time algorithm to test whether a graph contains an induced cycle with length at least $\ell$ and odd.
\end{abstract}

\section{Introduction}
All graphs in this paper are finite and have no loops or parallel edges. A {\em hole} of $G$ is an induced subgraph of $G$
that is a cycle
of length at least four, and an {\em antihole} is an induced subgraph whose complement is a cycle of length at least four. 
In 2005, two of us, with Cornu\'ejols, Liu and  Vu\v{s}kovi\'{c}~\cite{bergealg},
gave an algorithm to test whether an input graph $G$ has an odd hole or odd antihole, and thereby to 
test whether $G$ is perfect, with running time at most polynomial in $|G|$. ($|G|$ denotes the number of vertices of $G$.) 
At that time we were unable to separate the test for odd holes from the test for odd antiholes,
and testing for odd holes in poly-time has remained open until very recently. Indeed, it seemed quite likely that testing for an odd hole was NP-complete; for instance,
D. Bienstock~\cite{bienstock,bienstock2} showed that testing if a graph has an odd hole containing a given vertex is NP-complete.
So it was something of a surprise when recently we found a poly-time algorithm to test for odd holes~\cite{oddholetest}. 
(This is modified to run faster in a recent paper~\cite{thorup} by Lai, Lu, and Thorup.)

In this paper we extend that result: for each integer $\ell\ge 5$ we give a poly-time
algorithm to test whether $G$ has an odd hole of length at least~$\ell$. More exactly:

\begin{thm}\label{mainthm}
For each integer $\ell\ge 5$, there is an algorithm with the following specifications:
\begin{description}
\item [Input:] A graph $G$.
\item [Output:] Decides whether $G$ has an odd hole of length at least $\ell$.
\item [Running time:] $O(|G|^{20\ell+40})$.
\end{description}
\end{thm}
We have not tried very hard to optimize the exponent in the running time (although getting the exponent to be linear in $\ell$
took some effort).

We are not aware of previous work on detecting ``long'' induced subgraphs of specific type, although it seems a sensible question.
Here are three current pieces of work that are related:
\begin{itemize}
\item Linda Cook, with Seymour, has a poly-time algorithm to test if a graph has a long even hole~\cite{longevenhole}.
\item Eli Berger and Sophie Spirkl, with Seymour, have a poly-time algorithm to test if there is an induced path between 
specified vertices $s,t$ of a graph that has length longer than the shortest $st$-path~\cite{shortestpath}. It is open
whether there is a poly-time algorithm to test for an induced $st$-path of length at least three more than the shortest $st$-path.
\item We have a poly-time algorithm to find the shortest odd hole in a graph, if it has one~\cite{shortestoddhole}.
\end{itemize}
Long odd holes have been worked on before, although not for algorithms. In~\cite{oddholes} Scott and Seymour proved
that graphs with no odd hole are ``$\chi$-bounded '', and later, with Sophie Spirkl~\cite{longoddholes}, we extended this to graphs
with no long odd hole. Indeed, in~\cite{holeparity} Scott and Seymour extended it further, to graphs with no holes of length
$p$ modulo $q$, for any fixed $p,q$. We currently see no prospect of extending the algorithmic work to test for a hole of length
$p$ modulo $q$; even testing for holes of length a multiple of three seems very challenging, although of interest because such graphs
have nice properties~\cite{ternary}.

The new algorithm once again uses ``cleaning'', as does the algorithm of~\cite{oddholetest} and several other algorithms to detect special induced 
subgraphs. Indeed it was modelled on the algorithm of~\cite{oddholetest}, but it is considerably more complicated.

Here is an outline of the method. Throughout the paper, $\ell\ge 5$ is a fixed number, and throughout, a {\em long}  path or hole 
means a path or hole of length at least $\ell$. 
If $C$ is a hole in $G$, a vertex $v$ of $V(G)\setminus V(C)$ is {\em $C$-major} if there is no three-vertex path of $C$ containing all the neighbours
of $v$ in $V(C)$. A hole $C$ is {\em clean} if no vertex is $C$-major.
\begin{itemize}
\item First we test for the presence in the input graph $G$ of certain kinds of induced subgraphs (``short'' long odd holes, ``long pyramids'' and ``long jewels'')
that we can test for in polynomial time, and whose
presence would imply that $G$ contains a long odd hole. We call these three kinds of subgraphs ``easily-detected configurations''.
We may assume these tests are unsuccessful.
\item Second, we generate a ``cleaning list'', a list of polynomially-many subsets of $V(G)$, such that if $G$ has a long odd hole, and $C$
is a long odd hole of minimum length (a {\em shortest long odd hole}) then some set $X$ in the list contains all the $C$-major vertices and contains
no vertex of $C$ itself.  This relies on the fact that $G$ contains none of the easily-detected configurations.

\item Third, for each $X$ in the cleaning list, we test whether $G\setminus X$ has a clean shortest long odd hole. (More exactly, we give an algorithm that 
either decides that $G\setminus X$ has a long odd hole, or decides that $G\setminus X$ has no clean shortest long odd hole.)
This again relies on the absence of the easily-detected configurations.
\end{itemize}
The reader familiar with the method of~\cite{oddholetest} will see the similarity of the two algorithms.

But part of the approach is significantly different. To generate the cleaning list in \cite{oddholetest},
we used a theorem that if $C$ is a shortest odd hole, and $M$ is a set of $C$-major vertices such that one of them is 
nonadjacent to all the others, then there is a ``heavy edge'' in $C$, an edge $uv$ of $C$ such that every vertex in $M$ 
is adjacent to one of $u,v$. We tried to extend this to the long odd hole situation, but failed. For our purposes,
a ``heavy path'' of $C$ of bounded length (that is, such that every vertex in $M$ has a neighbour in the path) would be just
as good as a heavy edge;
and this extension might be true, but we were unable to prove it. In its
place we had to use a considerably more complicated method, proving that there is a bounded set of paths of $C$, 
each of bounded length, such that every vertex in $M$ has a neighbour in one of the paths; and we could only prove this
when the exceptional vertex of $M$ was carefully chosen.

The paper is organized as follows. First, we explain how to test for the easily-detected configurations; this is a straightforward adaptation  of the 
algorithms in~\cite{bergealg} to test for pyramids and jewels. Then we give the algorithm for the third step above; and finally we show how to 
generate the cleaning list.

Let us remark, finally, that if we want to test for a long hole, rather than a long odd hole, then this is easy: enumerate
all induced paths of $G$ of length $\ell-1$, and for each one, test directly if it can be extended to a hole. 
This has running time $O(|G|^{\ell+1})$, and so both this and our algorithm for \ref{mainthm} have running time
$|G|^{O(\ell)}$. We do not know if either can be substantially improved, although both problems are NP-hard when $\ell$ is part of 
the input. (To see this for the long hole problem, take a graph with $n$ vertices and subdivide each edge once, and take $\ell=2n$; 
then this has an induced cycle of length at least $\ell$
if and only if the original graph has a Hamilton cycle. For the long odd hole problem, take a graph with an odd number $n$ of vertices,
subdivide every edge twice, and take $\ell=3n$.)

\section{The easily-detected configurations}

We begin with the test for what we earlier called ``short'' long odd holes:

\begin{thm}\label{testshortlongoddhole}
There is an algorithm with the following specifications:
\begin{description}
\item [Input:] A graph $G$, and an integer $k\ge 0$.
\item [Output:] Decides whether there is a long odd hole in $G$ of length at most $k$.
\item [Running time:] $O(|G|^{k})$.
\end{description}
\end{thm}
\Proof We enumerate all sets of at most $k$ vertices, and for each one, check whether it induces a long odd hole.~\bbox

If $X\subseteq V(G)$, we denote the subgraph of $G$ induced on $X$ by $G[X]$. If $X$ is a vertex or edge of $G$, or a set of vertices
or a set of edges of $G$, we denote by $G\setminus X$ the graph obtained from $G$ by deleting $X$.
Thus, for instance, if $b_1b_2$ is an edge of a hole $C$, then $C\setminus \{b_1,b_2\}$ and $C\setminus b_1b_2$ are both paths,
but one contains $b_1,b_2$ and the other does not. If $P$ is a path,
we denote by $P^*$ the {\em interior} of $P$, the set of vertices of the path $P$ that are not ends of $P$.
If $P$ is a path and $x,y\in V(P)$,
we denote the subpath with ends $x,y$ by $x\dd P\dd y$. The {\em length} of a path or cycle is the number of edges in it.

Let $u,v\in V(G)$, and let $Q_1,Q_2$ be induced paths between $u,v$, of different parity.
Let $P$ be an induced path between $u,v$ of length at least $\ell$, such that no vertex of $P^*$ equals or is adjacent 
to any vertex of $Q_1^*,Q_2^*$.
We say the subgraph induced on $V(P\cup Q_1\cup Q_2)$ is a {\em long jewel} of {\em order} $\max(|V(Q_1)|,|V(Q_2)|)$, 
{\em formed by} 
$Q_1,Q_2,P$. Any graph containing a long jewel has a long odd hole, since the holes $P\cup Q_1,P\cup Q_2$ are both long and one of them is odd.
The next result extends theorem 3.1 of~\cite{bergealg}:

\begin{thm}\label{testlongjewel}
There is an algorithm with the following specifications:
\begin{description}
\item [Input:] A graph $G$, and an integer $k\ge 0$.
\item [Output:] Decides whether there is a long jewel in $G$ of order at most $k$.
\item [Running time:] $O(|G|^{2k+\ell})$.
\end{description}
\end{thm}
\Proof 
We enumerate all triples of induced paths $Q_1,Q_2,R$ of $G$, such that:
\begin{itemize}
\item $Q_1,Q_2$ join the same pair of vertices, say $u,v$;
\item one of $Q_1,Q_2$ is odd and the other is even, and each has at most $k$ vertices;
\item $R$ has length $\ell-2$, and has one end $u$ and the other some vertex $w$ say;
\item no vertex of $V(R)\setminus \{u\}$ equals or has a neighbour in $V(Q_1\cup Q_2)\setminus \{u\}$.
\end{itemize}
For each such triple of paths, let $X$ be the set of vertices of $G$ that are different from and nonadjacent to each vertex of 
$V(Q_1\cup Q_2\cup R)\setminus \{v,w\}$. We test whether there is a path in $G[X\cup\{w,v\}]$ between $w,v$. If so we output that $G$ contains
a long jewel of order at most $k$. If no triple yields this outcome, we output that $G$ has no such long jewel. 

To see the correctness of the algorithm, certainly the output is correct if $G$ contains no long jewel of order at most $k$. Suppose then it does, say 
formed by $Q_1,Q_2,P$. Let $u,v$ be the ends of $P$, and let $R$ be the subpath of $P$ of length $\ell-2$ with one end $u$. When the algorithm
tests the triple $Q_1,Q_2,R$, it will discover there is a path in $G[X\cup \{w,v\}]$ between $w,v$, because the remainder of $P$ is such a path. Consequently
the output is correct. 

The running time is $O(|G|^2)$ for each triple of paths, and there are at most $|G|^{2k+\ell-2}$ such triples, so the running time is as claimed.
This proves~\ref{testlongjewel}.~\bbox

Many of the algorithms in this paper follow the same outline; we enumerate all subgraphs, or sequences of vertices, of some prescribed type,
and for each one, perform some test on it. (Critically, there must be only polynomially many such subgraphs to test.)
If the test is successful, we have found a subgraph of the desired type, and if it is never successful
we will apply a theorem that says that then there is no subgraph of the desired type. 
For brevity we call the process of enumerating all these subgraphs and testing them one-by-one ``guessing'';
thus we would describe the long jewel algorithm above as ``guessing the two paths $Q_1,Q_2$ and an initial subpath of $P$''.

Let $v_0\in V(G)$, and for $i = 1,2,3$ let $P_i$ be an induced path of $G$ between $v_0$ and $v_i$, such that

\begin{itemize}
\item $P_1,P_2,P_3$ are pairwise vertex-disjoint except for $v_0$;
\item $v_1,v_2,v_3\ne v_0$, and at least two of $P_1,P_2,P_3$ have length at least $\ell$;
\item $v_1,v_2,v_3$ are pairwise adjacent; and
\item for $1\le i<j\le 3$, the only edge between $V(P_i)\setminus \{v_0\}$ and $V(P_j)\setminus \{v_0\}$ is the edge $v_iv_j$.
\end{itemize}
We call the subgraph induced on $V(P_1\cup P_2\cup P_3)$ a {\em long pyramid}, with {\em apex} $v_0$ and {\em base} $\{v_1,v_2,v_3\}$, {\em formed by}
$P_1,P_2,P_3$.
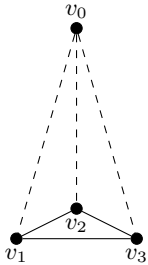
\begin{figure}[H]
\centering

\begin{tikzpicture}[scale=0.8,auto=left]
\tikzstyle{every node}=[inner sep=1.5pt, fill=black,circle,draw]

\node (v0) at (0,0) {};
\node (v1) at (-1,-3.5) {};
\node (v2) at (0,-3) {};
\node (v3) at (1,-3.5) {};
\tikzstyle{every node}=[]
\draw (v1) node [below]           {\footnotesize$v_1$};
\draw (v2) node [below]           {\footnotesize$v_2$};
\draw (v3) node [below]           {\footnotesize$v_3$};
\draw (v0) node [above]           {\footnotesize$v_0$};

\draw (v1) -- (v2);
\draw (v1) -- (v3);
\draw (v2) -- (v3);
\draw[dashed] (v0) -- (v1);
\draw[dashed] (v0) -- (v2);
\draw[dashed] (v0) -- (v3);

\end{tikzpicture}

\caption{A long pyramid. The dashed lines represent paths, of indeterminate length, but two of them must have length at least $\ell$.} \label{fig:pyramid}
\end{figure}
If $G$ has a long pyramid then $G$ has a long odd hole (because two of the paths $P_1,P_2,P_3$ have the same length modulo two, and they induce
a long odd hole). The next result extends theorem 2.2 of~\cite{bergealg}, and is proved similarly:

\begin{thm}\label{testlongpyramid}
There is an algorithm with the following specifications:
\begin{description}
\item [Input:] A graph $G$.
\item [Output:] Decides whether there is a long pyramid in $G$.
\item [Running time:] $O(|G|^{6\ell+8})$.
\end{description}
\end{thm}
\Proof
Suppose that $G$ contains a long pyramid; then it contains a ``smallest'' one, one with the fewest vertices, say with apex $v_0$ and 
base $\{v_1,v_2,v_3\}$, formed by
the path $P_1,P_2,P_3$. For $i = 1,2,3$, let $m_i$ be a vertex of $P_i$ that divides it into two paths with lengths differing by at most one.
For $i = 1,2,3$,
\begin{itemize}
\item if $P_i$ has length at least $\ell$, let $A_i$ be a subpath of $P_i$ of length $\ell$ with one end $v_0$, and let $B_i$
be a subpath of $P_i$ of length $\ell$ with one end $v_i$;
\item if $P_i$ has length less than $\ell$, let $A_i=B_i=P_i$.
\end{itemize}
Let $A_i$ have ends $v_0, a_i$, and let $B_i$ have ends $v_i,b_i$ for $i = 1,2,3$.

The algorithm proceeds as follows. If there is a pyramid as above, we guess the vertex $v_0$ and for $i=1,2,3$ we guess the vertices
$v_i,m_i,a_i,b_i$ and the paths $A_i,B_i$. Let $X$ be the set of all these vertices (including the vertices of the paths $A_i, B_i$ for $1\le i\le 3$).
For each $i$ such that $m_i$ does not belong to $V(A_i)$ we choose a shortest path $A_i'$ between $m_i,a_i$ such that its interior
is disjoint from $X\setminus \{m_i,a_i\}$ and 
contains no vertex with a neighbour in $X\setminus \{m_i,a_i\}$, and let $Q_i'=A_i\cup A_i'$. If $m_i\in V(A_i)$ let $Q_i'=A_i$.
Similarly, if $m_i\notin V(B_i)$ we choose a shortest path $B_i'$ between 
$m_i,b_i$ such that its interior
is disjoint from $X\setminus \{m_i,b_i\}$ and
contains no vertex with a neighbour in $X\setminus \{m_i,b_i\}$; and let $R_i'=B_i\cup B_i'$. If $m_i\in V(B_i)$ let $R_i'=B_i$.

Now for $1\le i\le 3$ we test whether $Q_i'\cup R_i'$ is an induced path between $v_0,v_1$, and if this is true for each $i$, and the three
paths form a pyramid, 
we return that $G$ contains a long pyramid. To prove the correctness, we must now prove a theorem that starting from a smallest pyramid
as described, $Q_i'\cup R_i'$ is indeed a path between $v_0,v_i$ for $1\le i\le 3$, and these three paths form a (possibly different) smallest pyramid.

Let $\Pi$ be a smallest pyramid in $G$, formed by paths $P_1,P_2,P_3$ as above.
For $i=1,2,3$, let $Q_i,R_i$ be the subpaths of $P_i$ between $m_i,v_0$ and between $m_i, v_i$ respectively. 
\\
\\
(1) {\em With $Q_1'$ chosen as in the algorithm, no vertex of $Q_1'$ belongs to or has a neighbour in $V(P_2\cup P_3)$.}
\\
\\
Suppose that this is false. Consequently $Q_1'\ne A_1$, and so $m_1\notin V(A_1)$ and $m_1$ is an end of $Q_1'$. 
Let $S$ be a minimal subpath of $Q_1'$ with one end $m_1$
such that its other end has a neighbour in $V(P_2\cup P_3)$. Let $S$ have ends $m_1,s$. There are three cases. 

First, suppose that $s$ has a unique neighbour $t$ in $P_2\cup P_3$.
We may assume that $t\in V(P_2)$ from the symmetry; let $\Pi'$ be the pyramid with apex $t$ and base $\{v_1,v_2,v_3\}$, 
formed by the paths $t\dd P_2\dd v_2$, $t\dd P_2\dd v_0\dd P_3\dd v_3$
and an induced path between $t,v_1$ with interior in $V(S\cup R_1)$. This is indeed a pyramid, from the minimality of $S$; it is long,
since all three of the paths have length at
least $\ell$ (because they include the paths $B_2, B_3,B_1$ respectively); and it has fewer vertices than $\Pi$, since $|V(S)|\le |V(Q_1)|-\ell$.
This is impossible from the choice of $\Pi$.

Second, suppose that $s$ has two nonadjacent neighbours in $V(P_2\cup P_3)$. 
Let $\Pi'$ be the pyramid with apex $s$ and base $\{v_1,v_2,v_3\}$
formed by the induced paths between $s,v_2$ and between $s,v_3$, both with interior in $V(P_2\cup P_3\}$ (these are unique, and include $B_2,B_3$ respectively, and so are long), and a path
between $s,v_1$ with interior in $V(S\cup R_1)$. Again, this is a long pyramid with fewer vertices than $\Pi$, a contradiction.

Third, suppose that $s$ has exactly two neighbours $t_1,t_2$ in $V(P_2\cup P_3)$ and they are adjacent. We may assume that 
$t_1,t_2\in V(P_2)$ and $t_2$ is closer to $v_2$ in $P_2$. Let $\Pi'$ be the pyramid with apex $v_0$ and base 
$\{s,t_1,t_2\}$ formed by the paths $v_0\dd P_2\dd t_1$, $v_0\dd P_3\dd v_3\dd v_2\dd P_2\dd t_2$ and a path
between $v_0,s$ with interior in $V(S\cup Q_1)$. Again, this is a long pyramid (because the three paths include $A_2,A_3,A_1$ respectively),
and it has fewer vertices than $\Pi$, because $|V(S)|\le |V(Q_1)|-\ell< |V(R_1)|$, a contradiction. This proves (1).
\\
\\
(2) {\em With $R_1'$ chosen as in the algorithm, no vertex of $R_1'$ belongs to or has a neighbour in $V(P_2\cup P_3)$.}
\\
\\
The proof is similar and we omit it.

\bigskip

Let $P_1'$ be an induced path between $v_0,v_1$ with interior in $V(Q_1'\cup R_1')$. 
From (1) and (2), it follows that 
$P_1', P_2,P_3$ form a long pyramid with at most as many vertices as $\Pi$. Consequently equality holds, and so 
$P_1'=Q_1'\cup R_1'$;
and $P_1',P_2,P_3$ form a smallest long pyramid. Similarly 
$Q_2'\cup R_2'$ is an induced path between $v_0,v_2$, say $P_2'$; and $P_1',P_2',P_3$ form a smallest long pyramid.
And similarly for $Q_3'\cup R_3'$. This proves the correctness of the algorithm. 

For its running time, we are guessing a sequence of at most $6(\ell+1)$ vertices, so the running time is as claimed. This proves \ref{testlongpyramid}.~\bbox

Let us say $G$ is a {\em candidate} if $G$ contains no long pyramid, no long jewel of order at most $\ell+2$, and no long odd hole
of length at most $2\ell+2$. In view of \ref{testshortlongoddhole}, \ref{testlongjewel}, and \ref{testlongpyramid}, we have:
\begin{thm}\label{testcandidate}
There is an algorithm with the following specifications:
\begin{description}
\item [Input:] A graph $G$.
\item [Output:] Decides whether $G$ is a candidate.
\item [Running time:] $O(|G|^{6\ell+8})$.
\end{description}
\end{thm}
Any graph that is not a candidate has a long odd hole, so now we just need to find a poly-time algorithm to test 
whether candidates have long odd holes.

\section{Detecting a clean shortest long odd hole}

The following was proved in~\cite{bergealg}:
\begin{thm}\label{oldshortpath}
Let $G$ be a graph containing no jewel or pyramid, and let $C$ be a clean shortest odd hole in $G$.
Let $u,v\in V(C)$ be distinct and nonadjacent, and let $L_1,L_2$ be the two subpaths of $C$ joining $u,v$,
where $|E(L_1)| < |E(L_2)|$. Then:
\begin{itemize}
\item  $L_1$ is a shortest path in $G$ between $u,v$, and
\item for every shortest path $P$ in $G$ between $u, v$, $P\cup L_2$ is a shortest odd hole in $G$.
\end{itemize}
\end{thm}
This was crucial in that paper. Happily, the exact analogue holds for clean shortest long odd holes:
\begin{thm}\label{shortpath}
Let $G$ be a graph containing no long jewel of order at most $\ell+2$, and no long pyramid, and with no long odd hole of 
length at most $2\ell+2$. Let 
$C$ be a clean shortest long odd hole in $G$.
Let $u,v\in V(C)$ be distinct and nonadjacent, and let $L_1,L_2$ be the two subpaths of $C$ joining $u,v$,
where $|E(L_1)| < |E(L_2)|$. Then:
\begin{itemize}
\item  $L_1$ is a shortest path in $G$ between $u,v$, and
\item for every shortest path $P$ in $G$ between $u, v$, $P\cup L_2$ is a clean shortest long odd hole in $G$.
\end{itemize} 
\end{thm}

First we prove the first assertion of \ref{shortpath}. 
If $u,v$ are vertices of a graph $G$, $d_G(u,v)$ denotes the length of the shortest path of $G$ joining $u,v$ ($d_G(u,v)=\infty$ if there is no such
path).
\begin{thm}\label{shortcut}
Let $G$ be a graph containing no long pyramid, no long jewel of order at most $\ell+2$, and no long odd hole of length at most $2\ell$. Let $C$
be a clean shortest long odd hole in $G$.
Then $d_G(u,v)=d_C(u,v)$ for all $u,v\in V(C)$.
\end{thm}
\Proof
Suppose the result is false; then there is an induced path $Q$ with vertices $q_1\cc q_k$ in order, such that some vertex $u\in V(C)$
is adjacent to $q_1$, some $v\in V(C)$ is adjacent to $q_k$, and $d_C(u,v)>k+1$. Choose such a path $Q$ with $k$ minimum.
Since $d_C(u,v)>k+1\ge 2$, and $q_1$ is not $C$-major, it follows that $q_1\ne q_k$, and so $k\ge 2$. If some vertex of $Q$ belongs to 
$V(C)$, say $q_i\in V(C)$, then from the choice of $k$, $d_C(u,q_i)\le i$, and $d_C(q_i,v)\le k-i+1$, and so
$$d_C(u,v)\le d_C(u,q_i)+d_C(q_i,v)\le k+1,$$ 
a contradiction. Thus $Q\cap C$ is null. 
There are two paths of $C$ that join $u,v$; one, $L_1$ say, of length $d_C(u,v)$, and the other, $L_2$ say,
 longer and of opposite parity.

Since $q_1$ is not $C$-major, there is a path $P_1$ of $C$ of length at most two such that all neighbours of $q_1$
in $V(C)$ lie in $V(P_1)$; choose $P_1$ minimal, and consequently both its ends are adjacent to $q_1$. (Possibly $P_1$
has only one vertex.) Define $P_2$
similarly for $q_k$.
\\
\\
(1) {\em $P_1, P_2$ are vertex-disjoint.}
\\
\\
Suppose that $P_1\cap P_2$ is non-null. Thus $d_C(u,v)\le 4$, and since $d_C(u,v)>k+1\ge 3$ it follows that $k=2$ and $d_C(u,v)=4$,
and $P_1,P_2$ both have length exactly two. Hence $P_1\cup P_2=L_1$, and the three paths $L_1$, $u\dd q_1\dd q_2\dd v$, and $L_2$,
form a long jewel of order four, a contradiction. Thus $P_1\cap P_2$ is null. This proves (1).

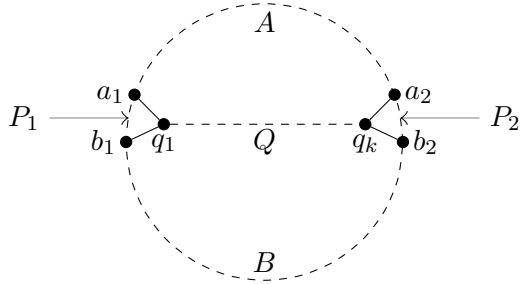
\begin{figure}[H]
\centering

\begin{tikzpicture}[scale=0.8,auto=left]
\tikzstyle{every node}=[inner sep=1.5pt, fill=black,circle,draw]

\def\r{2.3}
\def\s{1.7}
\node (b1) at ({-\r*cos(0)}, {\r*sin(0)}) {};
\node (b2) at ({-\r*cos(180)}, {\r*sin(180)}) {};
\node (a1) at ({-\r*cos(20)}, {\r*sin(20)}) {};
\node (a2) at ({-\r*cos(160)}, {\r*sin(160)}) {};
\node (q1) at ({-\s*cos(10)}, {\s*sin(10)}) {};
\node (qk) at ({-\s*cos(170)}, {\s*sin(170)}) {};
\draw (q1)--(a1);
\draw (q1)--(b1);
\draw(qk)--(a2);
\draw (qk) -- (b2);

\draw[domain=0:360,smooth,variable=\x,dashed] plot ({\r*cos(\x)},{\r*sin(\x)});

{\scriptsize}{
\draw[dashed] (q1)--(qk);
\tikzstyle{every node}=[]
\draw (b1) node [left]           {$b_1$};
\draw (a1) node [left]           {$a_1$};
\draw (b2) node [right]           {$b_2$};
\draw (a2) node [right]           {$a_2$};
\draw (q1) node [below]           {$q_1$};
\draw (qk) node [below]           {$q_k$};

\draw (0,0) node{$Q$};
\draw (0,2) node{$A$};
\draw (0,-2) node{$B$};

\node (P1) at (-4,{\r*sin(10)}) {$P_1$};
\draw[
    gray, ultra thin,decoration={markings,mark=at position 1 with {\arrow[black,scale=2]{>}}},
    postaction={decorate},
    ]
(P1) to ({-\r*cos(10)}, {\r*sin(10)}) {};
\node (P2) at (4,{\r*sin(170)}) {$P_2$};
\draw[
    gray, ultra thin,decoration={markings,mark=at position 1 with {\arrow[black,scale=2]{>}}},
    postaction={decorate},
    ]
(P2) to ({-\r*cos(170)}, {\r*sin(170)}) {};

}

\end{tikzpicture}

\caption{The paths $P_1,P_2,A,B$ of $C$.} \label{fig:shortcut}
\end{figure}

\bigskip
Let $P_1$ have ends $a_1,b_1$, and let $P_2$ have ends $a_2,b_2$,
where $a_1,b_1,b_2,a_2$ are in order in $C$ (possibly $a_1=b_1$ or $a_2=b_2$).
Let $A, B$ be two paths of $C$, where $A$ has
ends $a_1,a_2$, and $B$ has ends $b_1,b_2$, and the four paths $A,B,P_1,P_2$ are pairwise edge-disjoint and have union $C$.
\\
\\
(2) {\em One of $q_2\ll q_{k-1}$ has a neighbour in $V(C)$, and therefore $k\ge 3$.}
\\
\\
Suppose not. 
The hole formed by adding the path $a_1\dd q_1\cc q_k\dd a_2$ to $A$ is shorter than $C$, so either it has length less than $\ell$
or it has even length; so either $|E(A)|+k+1<\ell$, or $|E(A)|+k+1$ is even. Similarly either $|E(B)|+k+1<\ell$, or $|E(B)|+k+1$ is even. 
Since $|E(A)|+|E(B)|+4\ge |E(C)|>2\ell$, we may assume that $|E(B)|\ge \ell-1$.
In particular, $|E(B)|+k+1\ge \ell$, and so $|E(B)|+k$ is odd.
Now the paths $b_1\dd q_1\cc q_k\dd b_2$, $b_1\dd P_1\dd a_1\dd A\dd a_2\dd P_2\dd b_2$
and $B$ form a long jewel, of order the maximum of the lengths of the first two paths. But the second path has length
at least $d_C(u,v)\ge k+2$, so the second path is longer; and so the long jewel has order $|E(A)|+|E(P_1)|+|E(P_2)|$.
Consequently $|E(A)|+|E(P_1)|+|E(P_2)|>\ell+2$, and so $|E(A)|\ge \ell-1 $. Hence $|E(A)|+k+1\ge \ell$, and so  $|E(A)|+k$ is odd. Since
$C$ is odd, it follows that one of $P_1,P_2$ is odd and the other is even, and we may assume that $P_2$ has length one;
and so $q_k$ has exactly two neighbours in $V(C)$, $a_2$ and $b_2$, and they are adjacent. If $a_1=b_1$ then since $|E(A)|+|E(P_1)|+|E(P_2)|\ge \ell+3$, it follows
that $|E(A)|\ge \ell+2$, and so the three paths $A,B$ and $a_1\dd q_1\dd Q\dd q_k$ form a long pyramid, a contradiction.
Thus $P_2$ has exactly three vertices. Since $|E(A)|+|E(P_1)|+|E(P_2)|\ge \ell+3$, it follows that
$|E(A)|>\ell$, and so  the three paths $q_1\dd a_1\dd A\dd a_2$, $q_1\dd b_1\dd B\dd b_2$ and $Q$ form a long pyramid, a contradiction.
This proves (2).
\\
\\
(3) {\em None of $q_2\ll q_{k-1}$ has a neighbour in $V(L_2)$.}
\\
\\
Suppose that $q_i$ has a neighbour $w\in V(L_2)$ say, where $2\le i\le k-1$. 
Thus $w\ne u,v$. Let $R_2, S_2$ be the subpaths of $L_2$ between $w,u$ 
and between $w,v$ respectively. From the minimality of $k$, $d_C(u,w)<d_C(u,v)$, and so 
the path of $C$ between $u,w$ of length $d_C(u,w)$ is a subpath of $L_2$, that is, $R_2$ has length $d_C(u,w)$. 
Similarly $S_2$ has length $d_C(v,w)$, and so 
$d_C(u,w)+d_C(v,w)= |E(L_2)|$. 
From the minimality of $k$, $d_C(u,w)\le i+1$ (because otherwise $d_G(u,w)\le i+1<d_C(u,w)$ and $d_G(u,w)<k+1$, 
contrary to the minimality of $k$), and similarly
$d_C(v,w)\le k-i+2$, and so $|E(L_2)|\le k+3$. But $|E(L_2)|>|E(L_1)|=d_C(u,v)\ge k+2$, and so equality holds throughout.
In particular, $|E(L_2)|=|E(L_1)|+1=k+3$. Moreover, $R_2$ has length $d_C(u,w)=i+1$ and $S_2$ has length $k-i+2$.

The union of $L_1$ and the path $u\dd q_1\dd Q\dd q_k\dd v$
therefore is a cycle of length $|C|-2$. This is odd, less than $C$, and at least $\ell$, so this cycle is not induced.
Hence there exists $j\in \{1\ll k\}$ such that $q_j$ has a neighbour in the interior of $L_1$, say $x$. From the symmetry
under reversing $q_1\ll q_k$, 
we may assume that $j\ge i$. Let $R_1, S_1$ be the subpaths of $L_1$ between $x,u$ and between $x,v$ respectively.
The path $S_1\cup S_2$ has length more than the length of $S_2$
and hence at least $k-i+3$; and the path
$w\dd q_i\cc q_j\dd x$ has length $j-i+2\le k-i+2$. Thus $S_1\cup S_2$ is longer than $w\dd q_i\cc q_j\dd x$. But 
from the minimality of $k$, it follows that the length of $w\dd q_i\cc q_j\dd x$ is at least $d_C(w,x)$, and so at least the length
of $R_1\cup R_2$. Thus $j-i+2\ge |E(R_1)|+|E(R_2)|$. Also, the minimality of $k$ implies that the path $x\dd q_j\cc q_k\dd v$
has length at least the length of $S_1$, and so $k-j+2\ge |E(S_1)|$. Adding, we deduce that
$k-i+4\ge |E(R_1)|+|E(R_2)|+ |E(S_1)|$. But $|E(R_1)+|E(S_1)|=|E(L_1)|= k+2$, and $|E(R_2)|=i+1$, and so
$k-2i+4\ge k+3$, which is impossible. This proves (3).

\bigskip

We may assume that $u,v$ are chosen, adjacent to $q_1,q_k$ respectively, with $d_C(u,v)$ maximum. Then we have
\\
\\
(4) {\em $q_1,q_k$ have no neighbours in $L_2^*$.}
\\
\\
By (2) and (3), one of $q_2\ll q_{k-1}$ (say $q_i$) has a neighbour in $L_1^*$, say $w$. Since $d_G(u,w)<d_G(u,v)$, the minimality
of $k$ implies that $d_G(u,w)=d_C(u,w)$. If $w,v$ are adjacent, then 
$$d_C(u,v)\le d_C(u,w)+1=d_G(u,w)+1\le i+2\le k+1,$$
a contradiction; so $w$ is nonadjacent to $v$ and similarly nonadjacent to $u$. 
Consequently $w\in V(A)$.
Let $R,S$ be the subpaths of $A$ between $w$ and $a_1,a_2$ respectively.
Since $w, v$ are nonadjacent, the path $R\cup P_1$ has length at most that of $L_1$, and so less than $|E(C)|/2$, 
and hence $d_C(w,b_1)=|E(R\cup P_1)|$.
From the minimality of $k$, the path $b_1\dd q_1\cc q_i\dd w$ has length at least the length of $R\cup P_1$; and similarly
$w\dd q_i\cc q_k\dd b_2$ has length at least that of $S\cup P_2$. Adding, we deduce that
$P_1\cup A \cup P_2$ has length at most $k+3$. But $|E(L_1)|\ge k+2$, and we suppose for a contradiction that not both $b_1,b_2$
belong to $V(L_1)$, and so the length of $L_1$ is strictly less than that of $P_1\cup A\cup P_2$.
Hence we have equality throughout; so $|E(L_1)|=k+2$, and $P_1\cup A\cup P_2$ has length $k+3$, and 
exactly one edge of $P_1\cup P_2$ does not belong to $L_1$; and so
from the symmetry between $u,v$ we may assume that $u=b_1$ and $v,b_2$ are adjacent. But the path $b_1\dd q_1\cc q_k\dd b_2$ has length
$k+1$, and this has the same parity as, and is shorter by two than, the path $P_1\cup A\cup P_2$, so its union with the path $B$
is an odd hole $C'$ of length two less than that of $C$. Consequently $C'$ is not long; but this is impossible, 
since $C$ has length more than $2\ell$. This shows that $b_1,b_2\in V(L_1)$, and so proves (4).
\\
\\
(5) {\em $L_1$ has length $k+2$.}
\\
\\
From (2) and (3), there exists $i\in \{2\ll k-1\}$ such that $q_i$ has a neighbour $w\in L_1^*$. From the minimality of $k$,
the path $u\dd q_1\cc q_i\dd w$ has length at least the length of the subpath of $L_1$ between $u,w$; and 
the path $w\dd q_i\cc q_k\dd v$ has length at least the length of the subpath of $L_1$ between $w,v$. Adding, we deduce
that $k+3$ is at least the length of $L_1$. But adding the path $u\dd q_1\cc q_k\dd v$ to $L_2$ makes a hole, of length
at least $\ell$ since $L_2$ has length at least $|C|/2$; and this hole is shorter than $C$, and so has even length.
Consequently $k+1$ has the same parity as the length of $L_2$, and hence $k$ has the same parity as the length of $L_1$.
Since the length of $L_1$ is at most $k+3$, it equals $k+2$. This proves (5).

\bigskip

Let the vertices of $L_1$ in order be $u=w_0\dd w_1\dd w_{k+1}\dd w_{k+2}=v$.
\\
\\
(6) {\em For $1\le i\le k$ and $1\le j\le k+1$, if $q_i,w_j$ are adjacent then $j\in \{i,i+1\}$.}
\\
\\
If $j<i$ then the path $w_j\dd q_i\cc q_k\dd v$ has length less than the length of $w_j\dd w_{j+1}\cc w_{k+1}\dd v$, 
since the first path has length $k-i+2$ and the second has length $k+2-j=d_C(w_j,v)$. This is
contrary to the minimality of $k$. Similarly, if $j>i+1$, the path $u\dd q_1\cc q_i\dd w_j$ is shorter than $u\dd w_1\cc w_j$,
again a contradiction. This proves (6).

From (2), we may assume (exchanging $u,v$ if necessary) that there exists $i\in \{2\ll k-1\}$ such that $q_i$ is adjacent to $w_i$.
Let $C'$ be the union of $L_2$ and the path
$$u\dd w_1\cc w_i\dd q_i\cc q_k\dd v.$$
From (3), (4) and (6), $C'$ is a hole of length $|C|$, and so is a shortest long odd hole. Let $G'$ be the subgraph of $G$
induced on $V(C)\cup V(Q)$. Then $C'$ is a hole in $G'$, and moreover, in $G'$ $C'$ is a clean hole, because of (6).
But the path $u\dd q_1\cc q_i$ is shorter then the path $u\dd w_1\cc w_i\dd q_i$, and the latter has length $d_{C'}(u,q_i)$;
and this contradicts the minimality of $k$. This proves \ref{shortcut}.~\bbox

We need the following lemma.
\begin{thm}\label{minorvert}
Let $G$ be a graph containing no long jewel of order at most $k$, and no long odd hole of length less than $k+\ell$.
Let $C$ be a shortest long odd hole, and let $v\in V(G)$ be $C$-major.
Then every path of $C$ that contains
all the neighbours of $v$ in $V(C)$ has length more than $k$.
\end{thm}
\Proof
Suppose that $P$ is a path of $C$, with ends $a,b$ say, containing all the neighbours of $v$ in $V(C)$, and $P$ has length at most $k$.
Since $v$ is $C$-major, it follows that $P$ has length at least three.
Let $Q$ be the other path of $C$ with ends $a,b$. Since by hypothesis, $C$ has length at least $k+\ell$, it follows that $Q$
has length at least $\ell$. Adding the path $a\dd v\dd b$ to $Q$ therefore gives a long hole, and it is shorter than $C$ since $P$ has length
at least three. Consequently this hole is not odd; so $Q$ is even and so $P$ is odd. But then the
three paths $P$, $a\dd v\dd b$ and $Q$ form a long jewel of order at most $k$, a contradiction. This proves \ref{minorvert}.~\bbox

Now we prove the second statement of \ref{shortpath}, in the following.
\begin{thm}\label{reroute}
Let $G$ be a candidate, and let $C$
be a clean shortest long odd hole in $G$.
Let $u,v\in V(C)$ be nonadjacent, and let $Q$ be a shortest path in $G$ joining $u,v$. 
Let $L_1,L_2$ be the paths of $C$ that join $u,v$, where $L_1$ is shorter than $L_2$. 
Then $L_2\cup Q$ is a clean shortest long odd hole.
\end{thm}
\Proof
Let $Q$ have vertices $u\dd q_1\cc q_k\dd v$ in order. We proceed by induction on $k$. By \ref{shortcut},
$L_1$ and $Q$ have the same length. If some vertex $q_i$ of $Q^*$ belongs to $L_2^*$, then by two applications of \ref{shortcut},
$d_C(u,q_i) \le i$ and $d_C(q_i,v)\le k-i+1$, so $L_2$ has length at most $k+1$, which is impossible since $L_1$ has length $k+1$ and $L_2$ is longer.
Thus $Q^*\cap L_2^*=\emptyset$, and so $L_2\cup Q$ is a cycle, with the same length as $C$. 
\\
\\
(1) {\em $L_2\cup Q$ is induced.}
\\
\\
Suppose it is not induced. Then there exist $i\in \{1\ll k\}$ and $w\in L_2^*$ such that $q_i, w$ are adjacent. 
Let $R,S$ be the subpaths of $L_2$
between $w$ and $u,v$ respectively. From \ref{shortcut}, $d_C(u,w)\le i+1\le |E(L_1)|$, and so $R$ has length $d_C(u,w)\le i+1$
(because the other path of $C$ joining $u,w$ includes $L_1$ and so is too long);
and $d_C(w,v)\le k-i+2\le |E(L_1)|$, so similarly $S$ has length $d_C(u,w)\le k-i+2$. Thus
$$|E(L_2)|= |E(R)|+|E(S)|=k+3-(i+1-d_C(u,w))-(k+i-2-d_C(v,w)).$$
In particular, $L_2$ has length at most $k+3$. Since $L_1$ has length $k+1$ and $L_1,L_2$ have opposite parity,
it follows that $L_2$ has length $k+2$, and so 
$(i+1-d_C(u,w))+(k+i-2-d_C(v,w))=1$; 
and from the symmetry between $u,v$, 
we may assume that $d_C(w,v)= k-i+2$ and $d_C(u,w)= i$. 
Since the path $w\dd q_i\cc q_k\dd v$ has the same length as
$S$, and the latter is a shortest path between $v,w$ by \ref{shortcut}, it follows that
$w\dd q_i\cc q_k\dd v$ is also a shortest path between $v,w$. Suppose that $i>1$; then from the inductive hypothesis,
the union of the path $L_1\cup R$ and $w\dd q_i\cc q_k\dd v$ is a clean shortest long odd hole $C'$ say. 
The two subpaths of $C'$ between $u,q_i$ have lengths $|E(L_1)|+k-i = 2k+1-i$ and $|E(R)|+1$. By \ref{shortcut}, one of
these paths has length at most the length of $u\dd q_1\cc q_i$, that is, at most $i$. But $2k+1-i>i$ since $i\le k$; and
$|E(R)|+1>i$ since $|E(R)|=|E(L_2)|-|E(S)| = k+2-(k-i+2)=i$, a contradiction. This proves (1).

\bigskip

Let $C'= L_2\cup Q$; then $C'$ is a shortest long odd hole. It only remains to check that $C'$ is clean.
Suppose not.
Then there is a $C'$-major vertex $x$. Since $x$ is not $C$-major, $x$ has a neighbour in the interior of $Q$.
Since $Q$ is a shortest path, there is a subpath $P_1$ of $Q$ of length at most two containing all neighbours of $x$
in $V(Q)$; choose $P_1$ minimal. Since $x$ is $C'$-major, it has a neighbour in the interior of $L_2$. Since $x$ is not $C$-major,
there is a path $P_2$ of $L_2$, of length at most two, containing all neighbours of $x$ in $V(C)$; choose $P_2$ minimal.
Let $P_1$ have ends $a_1,b_1$, where $u,a_1,b_1,v$ are in order in $Q$, and let $P_2$ have ends $a_2,b_2$, where $u,a_2,b_2,v$
are in order in $L_2$. Let $A$ be the path of $C'$ between $a_1,a_2$ that contains $u$, and let $B$ be the path of $C'$
between $b_1,b_2$ that contains $v$. By \ref{minorvert}, the path $P_1\cup A\cup P_2$ has length at least $\ell+3$, 
and so does the path $P_1\cup B\cup P_2$. In particular, $P_1,P_2$ are vertex-disjoint and do not contain $u$ or $v$.
Since $C'$ has length at least $2\ell+3$, one of $A,B$ has length at least $\ell$, say $A$. Hence the hole obtained
by adding $a_1\dd x\dd a_2$ to $A$ is long, and shorter than $C$, so even; and hence $A$ has even length. The path
$P_1\cup B\cup P_2$ therefore has odd length, and so this path, $a_1\dd x\dd a_2$ and $A$ form a long jewel, of order the length
of $P_1\cup B\cup P_2$; and so this path has length at least $\ell+3$. If $P_1$ has length two, let it have vertices $q_{i-1}\dd q_i\dd q_{i+1}$; then the path 
$$u\dd q_1\cc q_{i-1}\dd v\dd q_{i+1}\cc q_k\dd v$$
is a path between $u,v$ with the same length as $Q$, and violates (1). Thus $P_1$ has length at most one. Consequently
$B$ has length at least $\ell$; and so $B$ is even. Since $C'$ is odd, and $A,B$ are both even, it follows that exactly one of $P_1,P_2$
is odd. If $P_2$ is odd, that is, if $a_2,b_2$ are adjacent, then since $P_1$ is even and has length at most one,
it follows that $a_1=b_1$; and the three paths $A,B$ and $a_1\dd x$ form a long pyramid. Thus $P_2$ is even and $P_1$ is odd.
If $a_2=b_2$ then similarly the three paths $A,B$ and $a_2\dd x$ form a long pyramid; so $P_2$ has length two.
Since $x$ is not $C$-major, it has no neighbours in $L_1$. 
Since $P_1$ is odd, it follows that $k\ge 2$, so from the inductive hypothesis, the hole obtained from $C$
by replacing the middle vertex of $P_2$ by $x$ is a clean shortest long odd hole, violating (1).
This proves that $C'$ is clean, and so completes the inductive proof of \ref{reroute}.~\bbox

Now we can give the main result of this section. 
(In the following, we could bring the running time down to $O(|G|^{4})$, but there is no need.)

\begin{thm}\label{testcleanhole}
There is an algorithm with the following specifications:
\begin{description}
\item [Input:] A candidate $G$.
\item [Output:] Decides either that $G$ has a long hole, or that there is no clean shortest long odd hole in~$G$.
\item [Running time:] $O(|G|^{5})$.
\end{description}
\end{thm}
\Proof If $C$ is a shortest long odd hole in $G$, let $u,v,w\in V(C)$ be chosen such that each of $d_C(u,v),d_C(v,w),d_C(w,u)$ equals
either $\lfloor |C|/3 \rfloor$ or $\lceil |C|/3 \rceil$. Here is the algorithm: guess $u,v,w$, find a shortest path between each 
pair of them, and test whether these three paths make a long odd hole. If so, output that $G$ has an odd hole. After checking all triples,
if none has produced an odd hole, output that $G$ has no clean shortest long odd hole. It follows immediately from \ref{shortpath}
that the output is correct.~\bbox

\section{Covering by a short path}
Now we begin the third step of the main algorithm: cleaning a shortest long odd hole. In this section we prove some preliminary lemmas.
Let $C$ be a shortest long odd hole in $G$, and let $x$ be $C$-major. An {\em $x$-gap} means a path of $C$ of length at least two,
with both ends adjacent to $x$ and with no interior vertices adjacent to $x$; so if $P$ is an $x$-gap then adding $x$ gives a hole.
If $x,y$ are distinct nonadjacent $C$-major vertices in $G$, an {\em $(x,y)$-gap} means a path $P$ of $C$ such that $V(P)$ is the interior
of an induced path of $G$ between $x$ and $y$. If $x,y$ are adjacent $C$-major vertices, an {\em $(x,y)$-gap} means a path
$P$ of $C$ such that $V(P)$ is the interior
of an induced path of $G\setminus e$  between $x$ and $y$, where $e$ is the edge $xy$.
Two sets of vertices $X,Y$ are {\em anticomplete} if they are disjoint and there are no edges between them. We use the same term
for two subgraphs $P,Q$; thus we say $P$ is {\em anticomplete} to $Q$ if $V(P)$ is anticomplete to $V(Q)$.

\begin{thm}\label{paircover}
Let $C$ be a shortest long odd hole in a candidate $G$, and let $u,v$ be nonadjacent $C$-major vertices. Then there
is a $(u,v)$-gap of length less than $\ell/2-1$.
\end{thm}
\Proof
Let $A, B$ be the sets of neighbours of $u,v$ in $V(C)$ respectively.
We may assume that $u,v$ have no common neighbour in $V(C)$, because that would make an $(x,y)$-gap of length zero; 
and so there are an even number of $(u,v)$-gaps. Let us number them
$D_1\ll D_{2k}$ say, in order in $C$. We may assume that each $D_i$ has length at least $\ell/2-1$, and in particular, have length
at least two, since $\ell\ge 5$. For $1\le i\le 2k$, let $D_i$ have
ends $a_i\in A$ and $b_i\in B$. We say $i,j\in \{1\ll 2k\}$ are {\em consecutive} if either $|j-i|=1$ or $|j-i|=2k-1$.
\\
\\
(1) {\em For all $i,j\in \{1\ll 2k\}$, if $D_i, D_j$ have different parity then $i,j$ are consecutive.}
\\
\\
Suppose not; then we may assume that $D_1, D_i$ have opposite parity where $3\le i\le 2k-1$. Let $a_1'$ be the vertex of $C$ adjacent to $a_1$ that
is not in $D_1$, and define $b_1'$ similarly. Since $i\ne 2,2k$ it follows that $a_1,b_1\notin V(D_i)$.
If $a_1'\in V(D_i)$ then since $i,1$ are not consecutive, it follows that $\{a_1,a_1'\}$ includes the vertex set of a $(u,v)$-gap; but
this is impossible since all $(u,v)$-gaps have length at least two. So $a_1', b_1'\notin 
V(D_i)$, and adding $u,v$ to $D_1\cup D_i$ gives a long odd hole, shorter than $C$, a contradiction. This proves (1).
\\
\\
(2) {\em $D_1\ll D_{2k}$ all have the same parity.}
\\
\\
Suppose not; then $k\le 2$, and we may assume that $D_i$ has the same parity as $i$ for $1\le i\le 2k$. If say $a_1=a_2$, then
adding $v$ to $D_1\cup D_2$ gives a long odd hole shorter than $C$ (note that $b_1\ne b_2$ and they are nonadjacent, since $u$
has at least two neighbours); so $D_1\ll D_{2k}$ are pairwise vertex-disjoint. Since $D_1,D_2$ have opposite parity, they are not
anticomplete, and we may assume that $a_1$ is adjacent to one of $a_2,b_2$; and not to $b_2$ since all $(u,v)$-gaps have length at least two. So $a_1,a_2$ are adjacent. Since $u$ is $C$-major, it follows that $k\ge 2$, and so $k=2$. Similarly either $a_2,a_3$ are adjacent, or $b_2,b_3$ are adjacent, and the first is impossible since $a_2$ is adjacent to $a_1$ and has a neighbour in the interior of $D_2$.
So $b_2b_3, a_3a_4, b_4b_1$ are edges; but then $C$ has even length, a contradiction. This proves (2).

\bigskip

From (2), an even number of edges of $C$ belong to $(u,v)$-gaps. Let $F$ be the graph obtained from $C$ by deleting the edges and internal vertices
of every $(u,v)$-gap. Since $C$ is odd, it follows that $|E(F)|$ is odd, and so some component $P$ of $F$ is odd. Now the ends of $P$
belong to two consecutive $(u,v)$-gaps, $D_1,D_2$ say; and so we may assume that $P$ has ends $a_1,a_2$. But then $D_1\cup P\cup D_2$ is odd, and adding $v$
makes a long odd hole of length less than $C$, a contradiction. This proves \ref{paircover}.~\bbox

A similar proof yields:
\begin{thm}\label{paircover2}
Let $C$ be a shortest long odd hole in a candidate $G$, and let $u,v$ be nonadjacent $C$-major vertices. Suppose that all $(u,v)$-gaps
are odd; then for every long $(u,v)$-gap $Q$, there is a $(u,v)$-gap of length at most $\ell-5$, anticomplete
to $Q$.
\end{thm}
\Proof Since an even number of edges belong to 
$(u,v)$-gaps, it follows as in the proof of \ref{paircover} (and using the notation of that proof) 
that some component $P$ of $F$ is odd, and we may assume that $P$ has ends $a_1,a_2$.
Then $D_1\cup P\cup D_2$ is odd, and adding $v$
makes an odd hole of length less than that of $C$, which therefore has length less than $\ell$.
Since $P$ is odd, it follows that the sum of the lengths of $D_1,D_2$
is at most $\ell-4$, and so each has length at most $\ell-5$. Suppose that neither of them is anticomplete to $Q$.
Since $Q$ exists and $Q\ne D_1,D_2$, there are at least four $(u,v)$-gaps;
let $Q=D_i$ say where $3\le i\le 2k$. Since $D_i, D_2$ are not anticomplete, it follows that $i\le 4$; and similarly
$i\ge 2k-1$, and so $k=2$, and we may assume that $i = 3$ from the symmetry. Since $D_1,D_3$ are not anticomplete, it follows that
$b_1=b_4$, $a_3=a_4$, and $a_4,b_4$ are adjacent. Also since $D_2,D_3$ are not anticomplete, the vertices $b_2,b_3$ are either 
equal or adjacent; and since $C$ is odd and $D_1,D_2,D_3,D_4$ and $P$ are odd, it follows that $b_2=b_3$. Now the three 
paths $D_1\cup P\cup D_2$, $b_4\dd v\dd b_2$ and $D_3$ form a long jewel of order $|V(D_1\cup P\cup D_2)|<\ell$, a contradiction.
This proves \ref{paircover2}.~\bbox

The next result is crucial, and provides the machinery behind all the proof of \ref{mainthm} that is novel. When we apply this, we will have
a subpath $P$ of some hole $C$, and some $C$-major vertices with neighbours in $V(P)$; the sets $A_1\ll A_k$ in \ref{markerlemma}
will be the neighbour sets in $V(P)$ of these $C$-major vertices.

\begin{thm}\label{markerlemma}
Let $P$ be a path of odd length, with ends $p,p'$. Let $A_1\ll A_k\subseteq V(P)$ be nonempty, and let us say a subpath $Q$ of
$P$ is ``covering'' if $V(Q)\cap A_s\ne \emptyset$ for $1\le s\le k$. Suppose that the minimal covering subpath with one end $p$,
and the minimal covering subpath with one end $p'$, have the same parity.
Then there is an odd subpath $Q$ of $P$ such that $Q$ is covering, and  for some (possibly equal) $s,t\in \{1\ll k\}$, one end of $Q$
belongs to $A_s$, the other end belongs to $A_t$, and $A_s,A_t$ contain no other vertex of $Q$.
\end{thm}
\Proof
Let $P$ have vertices $p_0\cc p_n$ in order. Thus $n$ is odd. Choose $d\in \{0\ll n\}$ minimum such that 
$p_0\cc p_d$ is covering. 
For $d\le j\le n$, choose $i\le j$ maximum such that $p_i\cc p_j$ is covering and 
define $m(j)=j-i$. If $m(d)$ is odd, then the subpath of $P$ between $p_{d-m(d)}$ and $p_d$ satisfies the theorem; so
we assume that $m(d)$ is even. Also, from the hypothesis, $m(n)$ and $d$ have the same parity,
so $m(n)+n$ and $m(d)+d$ have different parity. Consequently there exists $i$ with $d+1\le i\le n$ such that $m(i)+i$ has different parity
from $m(i-1)+i-1$.
It follows that $m(i), m(i-1)$ have the same parity. Choose $h\le i-1$ maximum such that
$p_h\cc p_{i-1}$ is covering. Consequently $i-1-h=m(i-1)$.
From the maximality of $h$, there exists $t\in \{1\ll k\}$ such that $p_h\in A_t$, and $p_{h+1}\ll p_{i-1}\notin A_t$.
Since $m(i), m(i-1)$ have the same parity, it follows that
one of $p_{h+1}\ll p_{i-1},p_i$ belongs to $A_t$,
and so $p_i\in A_t$. If $i,h$ have opposite parity, the path $p_h\cc p_i$ satisfies the theorem, so we assume they have the same parity.
Choose $g\le i$ maximum such that $p_g\cc p_i$ is covering.
Thus $i-g=m(i)$. From the maximality of $g$, there exists $s\in \{1\ll k\}$ such that 
$p_g\in A_s$ and $ p_{g+1}\ll p_{i}\notin A_s$. Since $m(i-1),m(i)$ have the same parity, it follows that
$i-1-h$ and $i-g$ have the same parity, that is, $g-h$ is odd. But $i,h$ have the same parity, so $i-g$ is odd, and the path
$p_g\cc p_i$ satisfies the theorem. This proves \ref{markerlemma}.~\bbox

\section{Bases}

In this section we prepare to apply \ref{markerlemma} to generate a cleaning list.

\begin{thm}\label{constantgap}
Let $C$ be a shortest long odd hole in a candidate $G$, and let $u,v$ be nonadjacent $C$-major vertices. 
There is not both a long odd $(u,v)$-gap and a long even $(u,v)$-gap.
\end{thm}
\Proof
Let $P$ be a long odd $(u,v)$-gap, and let $Q$ be a long even $(u,v)$-gap. They are not anticomplete,
since otherwise adding $u,v$ to their union gives a long odd hole, shorter than $C$, a contradiction.
If they share a vertex, then their union is either a long odd $u$-gap or a long odd $v$-gap, a contradiction. 
So they are vertex-disjoint. Let $P$ have ends $p_1,p_2$, and let $Q$ have ends $q_1, q_2$, where $u$ is adjacent
to $p_1,q_1$ and $v$ to $p_2,q_2$. We may assume that one of $q_1,q_2$ is adjacent to $p_2$. If $q_1$ is adjacent to $p_2$
then $v\dd p_2\dd q_1\dd Q\dd q_2\dd v$ is a long odd hole shorter than $C$. So $q_2$ is adjacent to $p_2$. Let $R$ be the
path of $C$ joining $p_1,q_1$ that does not contain $p_2,q_2$. Thus $R$ has odd length, and so the three paths 
$R$, $p_1\dd u\dd q_1$ and $p_1\dd P\dd p_2\dd q_2\dd Q\dd q_1$ form a long jewel; and therefore this jewel has order at least
$\ell+3$. Consequently $R$ has at least $\ell+4$ vertices. Let $r_1,r_2$ be the neighbours of $p_1,q_1$ respectively in $R$.
Suppose that $v$ has no neighbours in $V(R)\setminus \{p_1,r_1,r_2,q_1\}$. Since $v$ is $C$-major, it is adjacent to at least
one of $r_1,r_2$; not to exactly one, since it would make a long pyramid with $C$, and not with both since then
$v\dd r_1\dd R\dd r_2\dd v$ is a long odd hole shorter than $C$; in each case a contradiction. So $v$ 
has a neighbour in $V(R)\setminus \{p_1,r_1,r_2,q_1\}$. Now suppose that $u$ has no neighbour in this set. If $u$ is adjacent to neither
of $r_1,r_2$ then $u\dd p_1\dd R\dd q_1\dd u$ is a long odd hole shorter than $C$; if $u$ is adjacent to both
$r_1,r_2$ then $u\dd r_1\dd R\dd r_2\dd u$ is a long odd hole shorter than $C$; and if $u$ is adjacent to exactly one of $r_1,r_2$,
it makes a long pyramid with $C$, in each case a contradiction. So $u$ also has a neighbour in $V(R)\setminus \{p_1,r_1,r_2,q_1\}$.
But then there is a path joining $u,v$ with interior in $V(R)\setminus \{p_1,r_1,r_2,q_1\}$; and this path, with 
$u\dd p_1\dd P\dd p_2$ and $u\dd q_1\dd Q\dd q_2$, forms a long pyramid, a contradiction. This proves \ref{constantgap}.~\bbox

Let $C$ be a shortest long odd hole in a candidate $G$.
If $C$ is not clean, there is a maximal path $D$ of $C$ such that some $C$-major vertex ($x$ say) is adjacent to its ends and not to 
any of its internal vertices. We call $(x,D)$ a {\em base} (for $C$ in $G$).
A base $(x,D)$ is {\em remote} if $D$ has length at least $2\ell$, and no $C$-major vertex different from $x$
has a neighbour $w\in V(C)$ such that $d_C(w,d_i)<\ell$ for some $i\in \{1,2\}$, where $d_1,d_2$ are the ends of $D$.
(If there is a $C$-major vertex different from $x$, the second condition just given implies the first.)
It is easy to ``arrange'' algorithmically that a base is remote: we just guess the two paths of $C$ of length $2\ell$
 with middle vertex $d_1,d_2$ respectively, and delete all vertices not in these paths with a neighbour in the interior of 
one of them, except $x$. (This is safe, because no vertex of $C$ will be deleted.) If any $C$-major vertex different from $x$ remains, 
then the base has become remote (and if they all have been deleted then we have won). So the theorems to come will often assume that
$(x,D)$ is a remote base. If $(x,D)$ is a remote base then $D$ is even, since adding $v$ to $D$ gives a long hole shorter than $C$.

In what follows, when we have a base $(x,D)$, we will always denote the ends of $D$ by $d_1,d_2$. If $v$ is a $C$-major vertex
nonadjacent to $x$, then $v$ has a neighbour in $V(D)$ from the maximality of $D$, and since $(x,D)$ is remote,
there are two $(x,v)$-gaps included in $D$, $D_1(v)$ and $D_2(v)$ say, where for $i = 1,2$, 
$D_i(v)$ has ends $d_i$ and $d_i(v)$ say. Both of them are long, since $(x,D)$ is remote, so they have the same parity by \ref{constantgap}.

We will use this notation throughout the paper without defining it again.

\begin{thm}\label{xcons}
Let $C$ be a shortest long odd hole $C$ in a candidate $G$, and let $(x,D)$ be a remote base. 
Then for every $C$-major vertex $v$ nonadjacent
to $x$, 
all $(x,v)$-gaps have the same parity.
\end{thm}
\Proof 
Let $D'$ be the path of $C$ different from $D$ with ends $d_1,d_2$. Since $D$ is even (because adding $x$ gives a long hole
shorter than $C$) it follows that $D'$ is odd. Thus the three paths $D'$, $d_1\dd x\dd d_2$ and $D$ form a long jewel, and so
$D'$ has length at least $\ell+3$. Consequently $x$ has a neighbour in $C$ that is different from and nonadjacent to $d_1,d_2$
(to see this, suppose not; then $x$ has at most four neighbours in $C$, and if three then it makes a long pyramid, and if two or four
then it makes a long odd hole shorter than $C$). Hence every $(x,v)$-gap different from $D_1(v),D_2(v)$ is
anticomplete to one of $D_1(v),D_2(v)$, and therefore has the same parity as $D_1(v),D_2(v)$. This proves \ref{xcons}.~\bbox

Let $(x,D)$ be a remote base for $C$. If $v$ is a $C$-major vertex nonadjacent to $x$, we say the {\em $x$-parity} of $v$ is the 
common parity of all  
$(x,v)$-gaps. 

\begin{thm}\label{opengap}
Let $C$ be a shortest long odd hole in a candidate $G$, and let $(x,D)$ be a remote base.
Let $v_1,v_2$ be nonadjacent $C$-major vertices, nonadjacent to $x$, and with the same $x$-parity, 
and let $Q$ be a long odd $(v_1,v_2)$-gap edge-disjoint from $D$. Then $x$ has no neighbour in $V(Q)$.
\end{thm}
\Proof
Let $Q$ have ends $q_1,q_2$, where $v_1q_1$ and $v_2q_2$ are edges.
Since $q_i$ is adjacent to $v_i$, it follows that $d_C(q_i, d_j)\ge \ell$ for $j = 1,2$; and so $Q$ is anticomplete to $D$.

Suppose $x$ has a neighbour in $V(Q)$. If $x$ has no neighbour in the interior of $Q$,  then $x$ is adjacent to one of $q_1,q_2$, so
one and therefore both of $v_1,v_2$ have even $x$-parity, and consequently $Q$ is not an $(x,v_1)$ or $(x,v_2)$-gap; and 
therefore $x$ is adjacent to both $q_1,q_2$,
and adding $x$ to $Q$ would give a long odd hole shorter than $C$, a contradiction. So $x$ has a neighbour in the interior of $Q$.
\begin{figure}[H]
\centering

\begin{tikzpicture}[scale=0.8,auto=left]
\tikzstyle{every node}=[inner sep=1.5pt, fill=black,circle,draw]

\def\r{3}
\def\s{2.3}
\node (x) at (0,0) {};
\draw[domain=0:360,smooth,variable=\x,dashed] plot ({\r*cos(\x)},{\r*sin(\x)});
\node (v1) at ({\s*cos(135)},{\s*sin(135)}) {};
\node (v2) at ({\s*cos(45)},{\s*sin(45)}) {};
\node (q1) at ({\r*cos(135)},{\r*sin(135)}) {};
\node (q2) at ({\r*cos(45)},{\r*sin(45)}) {};
\draw (v1)--(q1);
\draw (v2)--(q2);
\node (d1) at ({\r*cos(200)},{\r*sin(200)}) {};
\node (d2) at ({\r*cos(-20)},{\r*sin(-20)}) {};
\draw (x) -- (d1);
\draw (x)--(d2);

\node (d1v1) at ({\r*cos(220)},{\r*sin(220)}) {};
\node (d0) at ({\r*cos(235)},{\r*sin(235)}) {};
\node (d1v2) at ({\r*cos(255)},{\r*sin(255)}) {};
\draw (v1) -- (d1v1);
\draw (v1) -- (d0);
\draw (v2) -- (d1v2);

\node (x1) at ({\r*cos(110)},{\r*sin(110)}) {};
\node (x2) at ({\r*cos(70)},{\r*sin(70)}) {};
\draw (x) -- (x1);
\draw (x) -- (x2);

\tikzstyle{every node}=[]
\draw (v1) node [right]           {$v_1$};
\draw (v2) node [left]           {$v_2$};
\draw (q1) node [left]           {$q_1$};
\draw (q2) node [right]           {$q_2$};
\draw (d1) node [left]           {$d_1$};
\draw (d2) node [right]           {$d_2$};
\draw (x) node [below]           {$x$};
\draw (d1v1) node [left]           {$d_1(v_1)$};
\draw (d0) node [below]           {$d_0$};
\draw (d1v2) node [below]           {$d_1(v_2)$};

\node (Q) at (0,2.7) {$Q$};
\draw[
    gray, ultra thin,decoration={markings,mark=at position 1 with {\arrow[black,scale=2]{>}}},
    postaction={decorate},
    ]
(Q) to[bend right=20] ({2.7*cos(135)}, {2.7*sin(135)}) {};
\draw[
    gray, ultra thin,decoration={markings,mark=at position 1 with {\arrow[black,scale=2]{>}}},
    postaction={decorate},
    ]
(Q) to[bend left=20] ({2.7*cos(45)}, {2.7*sin(45)}) {};

\node (D) at (0,-2.7) {$D$};
\draw[
    gray, ultra thin,decoration={markings,mark=at position 1 with {\arrow[black,scale=2]{>}}},
    postaction={decorate},
    ]
(D) to[bend left] ({2.7*cos(200)}, {2.7*sin(200)}) {};
\draw[
    gray, ultra thin,decoration={markings,mark=at position 1 with {\arrow[black,scale=2]{>}}},
    postaction={decorate},
    ]
(D) to[bend right] ({2.7*cos(-20)}, {2.7*sin(-20)}) {};

\node (P2) at ({3.3*cos(57)},{3.3*sin(57)}) {$P_2$};
\node (P1) at ({3.3*cos(123)},{3.3*sin(123)}) {$P_1$};


\end{tikzpicture}

\caption{For \ref{opengap}.} \label{fig:opengap}
\end{figure}

Consequently there is an $(x,v_1)$-gap $P_1$ that is anticomplete to $D$ and $v_2$ has no neighbour in $V(P_1)$; and also
there is $P_2$ similarly. For $i = 1,2$, let $P_i^+$ be the path between $v_i,x$ with interior $V(P_i)$.
Now every $(v_1,v_2)$-gap included in $D$ is anticomplete to $Q$, and since $Q$ is odd and has length at least $\ell$,
it follows that every $(v_1,v_2)$-gap included in $D$ has odd length. In particular, $d_1(v_1)\ne d_1(v_2)$.
From the symmetry we may assume that $D_1(v_2)$ is a proper subpath of $D_1(v_1)$. Choose a vertex $d_0$ of the path $d_1(v_1)\dd D\dd d_1(v_2)$,
adjacent to $v_2$, such that the subpath $d_1(v_1)\dd D\dd d_0$ is minimal.  The latter is a $(v_1,v_2)$-gap, and so odd. On the
other hand, since $v_1,v_2$ have same $x$-parity, the paths $D_1(v_1), D_1(v_2)$ have the same parity, 
and so $d_1(v_1)\dd D\dd d_1(v_2)$ is even; and hence 
$d_1(v_2)\dd D\dd d_0$ is odd.
If $d_1(v_2),d_0$ are nonadjacent, then the path $d_0\dd D\dd d_1(v_1)\dd v_1\dd P_1^+\dd x\dd d_1\dd D\dd d_1(v_1)$ can be extended to a hole
by adding either the path $d_1(v_2)\dd D\dd d_0$ or $d_1(v_2)\dd v_2\dd d_0$, and since these two paths have different parity, one of these holes
is odd. It is long since $d_1\dd D\dd d_1(v_2)$ has length at least $\ell$; and it is shorter than $C$, a contradiction. Thus $d_1(v_2),d_0$
are adjacent. Choose a minimal subpath $R$ of $D$ with one end $d_2$ such that the other end, $r$ say, is adjacent to one of $v_1,v_2$.
Thus $R$ has length at least $\ell$, and only one of $v_1,v_2$ has a neighbour in $V(R)$. If $r$ is adjacent to $v_1$, the three paths
$x\dd d_2\dd R\dd r\dd v_1\dd d_1(v_1)\dd D\dd d_0$, $x\dd d_1\dd D\dd d_1(v_2)$ and $x\dd P_2^+\dd v_2$ form a long pyramid (omitting $v_1$
from the first if $d_1(v_1),d_0$ are equal or adjacent). If $r$ is adjacent to $v_2$, the three paths
$x\dd d_2\dd R\dd r\dd v_2$, $x\dd d_1\dd D\dd d_1(v_1)$ and $x\dd P_1^+\dd v_1\dd d_1(v_1)\dd D\dd d_0$ form a long pyramid, a contradiction.
This proves \ref{opengap}.~\bbox

\section{Catch and clean}

Let $C$ be a shortest long odd hole, and let $P$ be a path of $C$. If $v$ is a vertex, we say that $P$ {\em catches} $v$
if $v\notin V(P)$ and $v$ is adjacent to an internal vertex of $P$. The point is, if $P$ is a path of $C$ then we are sure 
that no vertices it catches belong to $V(C)$; they might not all be $C$-major, but it does no harm to delete them.
This is a quite effective way to clean $C$ of $C$-major vertices:
we guess a path of bounded length (or a bounded number of such paths) and delete all the vertices each one catches. For instance, if 
we could prove that there is such a set of paths of $C$ (of bounded size, and each of bounded length) 
that together catch all the $C$-major vertices, we could clean $C$ and be done. 
We have not been able
to prove or disprove this. Nonetheless, catching $C$-major vertices by paths is a useful technique, and we will develop it in this section.
If $\mathcal{F}$ is a set of paths of $C$, its {\em cost} is the number of vertices in the 
union of the paths; and it {\em catches} $v$ if one of its paths catches $v$.

\begin{thm}\label{cleanoutadj}
Let $C$ be a shortest long odd hole in a candidate $G$, and let $(x,D)$ be a remote base.
Let $M$ be a set of $C$-major vertices, all nonadjacent to $x$ and with the same $x$-parity.
Let $v_1,v_2\in M$ be adjacent, 
and let $Q$ be an odd $(v_1,v_2)$-gap, with ends $q_1,q_2$, edge-disjoint from $D$,
where
$v_i, q_i$ are adjacent for $i = 1,2$.
Suppose that $x$ has a neighbour in $V(Q)$ and every vertex in $M$ has a neighbour in $V(Q)$.
Then there is a set of paths of $C$ with cost at most $5\ell +12$
that catches all the vertices in $M$.
\end{thm}
\begin{figure}[H]
\centering

\begin{tikzpicture}[scale=0.8,auto=left]
\tikzstyle{every node}=[inner sep=1.5pt, fill=black,circle,draw]

\def\r{3}
\def\s{2.3}
\node (x) at (0,0) {};
\draw[domain=0:360,smooth,variable=\x,dashed] plot ({\r*cos(\x)},{\r*sin(\x)});
\node (v1) at ({\s*cos(135)},{\s*sin(135)}) {};
\node (v2) at ({\s*cos(45)},{\s*sin(45)}) {};
\node (q1) at ({\r*cos(135)},{\r*sin(135)}) {};
\node (q2) at ({\r*cos(45)},{\r*sin(45)}) {};
\draw (v1)--(q1);
\draw (v2)--(q2);
\draw (v1) -- (v2);
\node (d1) at ({\r*cos(200)},{\r*sin(200)}) {};
\node (d2) at ({\r*cos(-20)},{\r*sin(-20)}) {};
\draw (x) -- (d1);
\draw (x)--(d2);

\node (d1v1) at ({\r*cos(250)},{\r*sin(250)}) {};
\draw (v1) -- (d1v1);
\draw (v2) -- (d1v1);

\node (p1) at ({\r*cos(100)},{\r*sin(100)}) {};
\node (p2) at ({\r*cos(80)},{\r*sin(80)}) {};
\draw (x) -- (p1);
\draw (x) -- (p2);

\tikzstyle{every node}=[]
\draw (v1) node [below left]           {$v_1$};
\draw (v2) node [below right]           {$v_2$};
\draw (q1) node [left]           {$q_1$};
\draw (q2) node [right]           {$q_2$};
\draw (d1) node [left]           {$d_1$};
\draw (d2) node [right]           {$d_2$};
\draw (x) node [below]           {$x$};
\draw (d1v1) node [below left]           {$d_1(v_1)=d_1(v_2)$};
\draw (p1) node[above] {$p_1$};
\draw (p2) node[above] {$p_2$};

\node (Q) at (0,2.7) {$Q$};
\draw[
    gray, ultra thin,decoration={markings,mark=at position 1 with {\arrow[black,scale=2]{>}}},
    postaction={decorate},
    ]
(Q) to[bend right=20] ({2.7*cos(135)}, {2.7*sin(135)}) {};
\draw[
    gray, ultra thin,decoration={markings,mark=at position 1 with {\arrow[black,scale=2]{>}}},
    postaction={decorate},
    ]
(Q) to[bend left=20] ({2.7*cos(45)}, {2.7*sin(45)}) {};

\node (D) at (0,-2.7) {$D$};
\draw[
    gray, ultra thin,decoration={markings,mark=at position 1 with {\arrow[black,scale=2]{>}}},
    postaction={decorate},
    ]
(D) to[bend left] ({2.7*cos(200)}, {2.7*sin(200)}) {};
\draw[
    gray, ultra thin,decoration={markings,mark=at position 1 with {\arrow[black,scale=2]{>}}},
    postaction={decorate},
    ]
(D) to[bend right] ({2.7*cos(-20)}, {2.7*sin(-20)}) {};

\node (P2) at ({3.3*cos(65)},{3.3*sin(65)}) {$P_2$};
\node (P1) at ({3.3*cos(115)},{3.3*sin(115)}) {$P_1$};


\end{tikzpicture}

\caption{For \ref{cleanoutadj}.} \label{fig:cleanoutadj}
\end{figure}
\Proof
For $i = 1,2$, let $P_i$ be the $(x,v_i)$-gap included in $Q$, and let its ends be $q_i, p_i$.
\\
\\
(1) {\em $d_1(v_1)=d_1(v_2)$.}
\\
\\
Suppose not; then we can assume that $D_1(v_1)$ is a proper subpath of $D_1(v_2)$. But then the holes
$$v_2\dd d_1(v_2)\dd D_1(v_2)\dd d_1\dd x\dd p_2\dd P_2\dd q_2\dd v_2$$ 
$$v_2\dd v_1\dd d_1(v_1)\dd D_1(v_1)\dd d_1\dd x\dd p_2\dd P_2\dd q_2\dd v_2$$
have opposite parity, and are both long and shorter than $C$, a contradiction. This proves (1).
\\
\\
(2) {\em $p_1,p_2$ are distinct; and if $p_1,p_2$ are nonadjacent then the sum of the lengths of $P_1,P_2$ is at most $\ell-6$.}
\\
\\
Since $P_1,P_2$ have the same parity, it follows that the path $p_1\dd Q\dd p_2$ is odd, and in particular $p_1\ne p_2$.
If $p_1,p_2$ are nonadjacent, there is an odd hole 
$$x\dd p_1\dd P_1\dd q_1\dd v_1\dd v_2\dd q_2\dd P_2\dd p_2\dd x$$
which is shorter than $C$, and hence has length less than $\ell$. This proves (2).
\\
\\
(3) {\em Let $T$ be a $(v_1,x)$-gap. If $v\in M$ is adjacent to $v_1$, then either $d_1(v)=d_1(v_1)$, or $v$ has a neighbour in $V(T)$.}
\\
\\
Let $v\in M$ be adjacent to $v_1$, and suppose that $d_1(v)\ne d_1(v_1)$, and $v$ has no neighbour in $V(T)$.
One of $D_1(v), D_1(v_1)$ is a proper subpath of the other. If $D_1(v_1)$ is a proper subpath of $D_1(v)$, choose a path $S$
joining $v,x$ with interior in the interior of $Q$; then since $S$ and $D_1(v_1)$ have the same parity, it follows that
the hole
$v\dd v_1\dd d_1(v_1)\dd D_1(v_1)\dd d_1\dd x\dd S\dd v$ is a long odd hole shorter than $C$, a contradiction.
If $D_1(v)$ is a proper subpath of $D_1(v_1)$, let $T^+$ be the path between $v_1,x$ with interior $V(T)$; then since $T$ and $D_1(v)$ have the same parity, the hole
$v_1\dd v\dd d_1(v)\dd D_1(v)\dd d_1\dd x\dd T^+\dd v_1$ is a long odd hole shorter than $C$, a contradiction. This proves (3).

\bigskip
For each $v\in M$ and for $i = 1,2$, let $Q_i(v)$ be the $(v,v_i)$-gap in $Q$, and let $q_i(v)$ be the neighbour of $v$ in $Q_i(v)$.
Now there are two cases, depending whether $p_1,p_2$ are adjacent or not.
\\
\\
(4) {\em If $p_1,p_2$ are adjacent, we can catch $M$ with a set of paths of $C$ with cost at most $5\ell +12$.}
\\
\\
Let $M_0$ be the set of all $v\in M$ such that either one of $Q_1(v), Q_2(v)$ has length less than $\ell$, or one of the $(x,v)$-gaps in $Q$ has length
less than $\ell$.
We can catch $M_0$ with a set of three paths of $C$, two of length $\ell+1$ and one of length $2\ell+1$, so with 
cost at most $4\ell+6$. For $i = 1,2$, let $M_i$ be the set of $v\in M\setminus M_0$ that are adjacent to $v_i$. By \ref{paircover}, there is a
$(v_1,x)$-gap $T$ of length less than $\ell/2-1$, and the same for $v_2$, so by (3) we can catch $M_1\cup M_2$ with a set of paths of $C$ with cost 
at most $\ell+6$. 
We claim that $M_0\cup M_1\cup M_2=M$; for suppose that $v\in M\setminus (M_0\cup M_1\cup M_2)$. 
Hence $v$ is nonadjacent to both $v_1,v_2$. Suppose first that $v$ has a neighbour in $P_1$
and a neighbour in $P_2$. Then there are two $(x,v)$-gaps in $Q$, with the same parity, and there is a long odd hole shorter than $C$, consisting of the union
of the two $(x,v)$-gaps, the edge $p_1p_2$, and the two edges from $v$ to the ends of the $(x,v)$-gaps; and this is impossible.
Thus we may assume that all neighbours of $v$ in $V(Q)$ belong to $V(P_1)$.
Since $x$ has neighbours in $Q_2(v)$, and $Q_2(v)$ is long, it is even by \ref{opengap}; but then
the $(v,x)$-gap and the $(x,v_2)$-gap in $Q_2(v)$ have different parities, a contradiction. This proves that  $M_0\cup M_1\cup M_2=M$. Consequently
we can catch $M$ with a set of paths of $C$ with cost at most $5\ell+12$. This proves (4).

\bigskip 

Henceforth we assume that $p_1,p_2$ are nonadjacent. 
Let $M_0$ be the set of all vertices $v\in M$ such that either one of $Q_1(v), Q_2(v)$ has length less than $\ell$, or $d_1(v)=d_1(v_1)$. 
We can catch $M_0$ with a set of three paths of $C$, with lengths $\ell+1,\ell+1$ and $2$, and hence with cost at most $2\ell+7$. 
\\
\\
(5) {\em If $v\in  M\setminus M_0$, then 
\begin{itemize}
\item $v$ is
nonadjacent to $v_1,v_2$
\item $Q_1(v), Q_2(v)$ are both even; and
\item $q_1(v), q_2(v)$ are adjacent (and so $v$ has exactly two neighbours in $V(Q)$, namely $q_1(v), q_2(v)$).
\end{itemize}
}
\noindent 
Since $Q_1(v)$ has length at least $\ell$,
and $P_1$ has length less than $\ell$ by (2), it follows that $v$ has no neighbour in $V(P_1)$. Since $P_1$ is an $(x,v_1)$-gap, and
$d_1(v)\ne d_1(v_1)$, (3) implies
that $v,v_1$ are nonadjacent, and similarly $v,v_2$ are nonadjacent. All neighbours of $v$ in $V(Q)$ lie strictly between $p_1,p_2$ in $Q$.
Since $Q_1(v)$ is long and contains a neighbour of $x$, \ref{opengap} implies that $Q_1(v)$ is even, and similarly $Q_2(v)$ is even. Thus $q_1(v)\ne q_2(v)$.
If $q_1(v), q_2(v)$ are nonadjacent, the hole 
$$v\dd q_1(v)\dd Q_1(v)\dd q_1\dd v_1\dd v_2\dd q_2\dd Q_2(v)\dd q_2(v)\dd v$$ 
is long, odd, and shorter
than $C$, a contradiction. Thus $q_1(v), q_2(v)$ are adjacent. This proves (5).
\\
\\
(6) {\em We can catch $M\setminus M_0$ with a set of paths of $C$ with cost at most $\ell+4$.}
\\
\\
We may assume that $M\setminus M_0$ is nonempty. Choose $u\in M\setminus M_0$ such that the path $Q_1(u)$ is as short as possible. Now there are three kinds
of vertices $v\in M\setminus M_0$:
\begin{itemize}
\item Let $N_0$ be the set of all $v\in M\setminus M_0$ such that some $(u,v)$-gap in $Q$ has length less than $\ell$ (the $(u,v)$-gap in $Q$ is unique unless
$u,v$ have the same two neighbours in $V(Q)$);
\item Let $N_1$ be the set of all $v\in M\setminus (M_0\cup N_0)$ such that $u,v$ are nonadjacent;
\item Let $N_2$ be the set of all $v\in M\setminus (M_0\cup N_0)$ such that $u,v$ are adjacent.
\end{itemize}
We can catch $N_0$ with one path of length $\ell+1$ (because of the minimality of $Q_1(u)$). We claim that $N_1=\emptyset$. 
For suppose that $v\in N_1$. Since $Q_1(u), Q_1(v)$ are both
even, it follows that the $(u,v)$-gap in $Q$ is odd, and long; and so $x$ has no neighbour in it, by \ref{opengap}. But then there is a long
pyramid formed by the path $x\dd d_1\dd D_1(v)\dd d_1(v)\dd v$, and the two $(x,v)$-gaps in $Q$, extended to $x$ (note that 
one of the latter is long, since
it includes the $(u,v)$-gap). Thus $N_1=\emptyset$. 

We claim that $d_1(v)=d_1(u)$ for all $v\in N_2$. Suppose not. If $D_1(u)$ is a proper subpath of $D_1(v)$, let $S$ be an induced path 
between $x, v$ with interior
in $V(Q)$ that is edge-disjoint from the $(u,v)$-gap, and therefore contains no neighbour of $u$. Since $D_1(u), D_1(v)$ have the same parity,
the holes 
$$v\dd d_1(v)\dd D_1(v)\dd d_1\dd x\dd S\dd v$$
$$v\dd u\dd d_1(u)\dd D_1(u)\dd d_1\dd x\dd S\dd v$$
have opposite parity, and both are long and shorter than $C$, a contradiction. If $D_1(v)$ is a proper subpath of $D_1(u)$, let $S$ be an induced path  
between $x, u$ with interior
in $V(Q)$ that is edge-disjoint from the $(u,v)$-gap, and therefore contains no neighbour of $v$. Since $D_1(u), D_1(v)$ have the same parity,
the holes
$$u\dd d_1(u)\dd D_1(u)\dd d_1\dd x\dd S\dd u$$
$$u\dd v\dd d_1(v)\dd D_1(v)\dd d_1\dd x\dd S\dd u$$
have opposite parity, and are both long and shorter than $C$, a contradiction. This proves that $d_1(v)=d_1(u)$ for all $v\in N_2$; and hence we can catch
$N_2$ with one path of length two and hence cost three. 
Since $N_1=\emptyset$ and we can catch $N_0$ with one path of length $\ell+1$ and  hence cost $\ell+2$, this proves (6).

\bigskip

But we can catch $M_0$ with a set of paths of $C$ with cost at most $2\ell+7$; so from (6) we can catch $M$ with cost at most $3\ell+11$.
This proves \ref{cleanoutadj}.~\bbox

As a consequence we have:
\begin{thm}\label{cleaneven}
Let $C$ be a shortest long odd hole in a candidate $G$, and let $(x,D)$ be a remote base. Then there is a set of paths of $C$ with
cost at most $5\ell+12$ that catches all $C$-major vertices that both are nonadjacent to $x$ and have even $x$-parity.
\end{thm}
\Proof
Let $M$ be the set of all $C$-major vertices nonadjacent to $x$ with even $x$-parity. Let $P$ be the path of $C$ 
obtained by deleting the interior of $D$. Every vertex in $M\cup \{x\}$ has at least two neighbours in $V(P)$. Moreover,
for each $v\in M\cup \{x\}$, the shortest subpath of $P$ with one end $d_1$ and the other adjacent to $v$ has
the same parity as $D_1(v)$ and so is even; and the same holds for $d_2$. Thus $M\cup \{x\}$ and $P$ satisfy the hypotheses
of \ref{markerlemma}, and so 
there is an odd path $Q$
of $C$ with ends $q_1,q_2$, edge-disjoint from $D$, and there are vertices $v_1,v_2\in M\cup \{x\}$, such that 
$v_1q_1$ and $v_2q_2$ are edges and there are no other edges between $\{v_1,v_2\}$ and $V(Q)$; and every vertex in $M\cup \{x\}$
has a neighbour in $V(Q)$. If $Q$ has length less than $\ell$ the result follows, so we assume $Q$ is long. Hence $v_1\ne v_2$.
Now $x$ is nonadjacent to all vertices in $M$, and they all have even $x$-parity; so $v_1,v_2\ne x$. By \ref{opengap}, $v_1,v_2$
are adjacent, and the result follows from \ref{cleanoutadj}. This proves \ref{cleaneven}.~\bbox

\begin{thm}\label{cleannonnbrs}
Let $C$ be a shortest long odd hole in a candidate $G$, and let $(x,D)$ be a remote base. Then there is a set of paths of $C$ with
cost at most $16\ell+28$ that catches all $C$-major vertices nonadjacent to $x$.
\end{thm}
\Proof
By \ref{cleaneven} there is a set $\mathcal{F}_1$ of paths of $C$ with
cost at most $5\ell+12$ that catches all $C$-major vertices nonadjacent to $x$ with even $x$-parity. Let $M_1$
be the set of all $C$-major non-neighbours of $x$ not caught by $\mathcal{F}_1$.
Now let us apply \ref{markerlemma} to the path $P$ of $C$
obtained by deleting the interior of $D$, and the set $M_1\cup \{x\}$. We may assume that $M_1\ne \emptyset$.
Again, the hypotheses of \ref{markerlemma} are satisfied,
since for each $v\in M_1$, the shortest subpath of $P$ with one end $d_1$ and the other adjacent to $v$ has
the same parity as $D_1(v)$ and so is odd; and the same holds for $d_2$; and so for $i = 1,2$ the shortest subpath of $P$ with 
one end $d_i$ containing a neighbour of each vertex in $M_1\cup \{x\}$ is odd. (Note that although $x$ is adjacent to $d_i$, 
the members of $M_1$
are not adjacent to $d_i$, and $M_1\ne \emptyset$, so this shortest subpath does not have length zero.)

It follows that there is an odd path $Q$
of $C$ with ends $q_1,q_2$, edge-disjoint from $D$, and there are vertices $v_1,v_2\in M\cup \{x\}$, such that
$v_1q_1$ and $v_2q_2$ are edges and there are no other edges between $\{v_1,v_2\}$ and $V(Q)$; and every vertex in $M_1\cup \{x\}$
has a neighbour in $V(Q)$. If $Q$ has length less than $\ell$ the result follows, so we assume $Q$ is long. Hence $v_1\ne v_2$.

First, suppose that $v_1,v_2\ne x$. 
By \ref{opengap}, $v_1,v_2$ are adjacent;
and by \ref{cleanoutadj} the result holds. So we may assume that $v_2=x$ say.
Let $B$ be the $v_1$-gap that includes $Q$. Then $(v_1,B)$ is a base for $C$ in $G\setminus N_1$,
where $N$ is the set of all $C$-major neighbours adjacent to $x$, and $N_1$ is the union of $N$ with the set of 
all $C$-major vertices caught by $\mathcal{F}_1$.
Let $\mathcal{F}_2$ be the set of the two paths of $C$ of length $2\ell$ with middle vertices the ends of $B$, and let $N_2$
be the union of $N_1$ with the set of $C$-major vertices caught by $\mathcal{F}_2$; then $(v_1,B)$ is a remote base for $C$ in $G\setminus N_2$.
By \ref{cleaneven} there is a set $\mathcal{F}_3$ of paths of $C$ with cost at most $5\ell+15$ that catches all 
$C$-major vertices of $G\setminus N_2$ nonadjacent to $v_1$ and with even $v_1$-parity.
Let $N_3$ be the union of $N_2$ with the set of $C$-major vertices caught by $\mathcal{F}_3$. Thus every $C$-major vertex in 
$G\setminus N_3$ is nonadjacent to $x$, and has odd $x$-parity, and odd $v_1$-parity, and its neighbours in $V(C)$ all
have $C$-distance
at least $\ell$ from each end of $B$.

By \ref{paircover2} there is an $(x,v_1)$-gap $T$ of length at most $\ell-5$ anticomplete to $Q$; and it follows that $T$
is vertex-disjoint from and anticomplete to $D$.
let $T^+$ be the path between $x,v_1$ with interior $V(T)$.
Let $T$ have ends $t_1,t_2$, where $t_1v$ and $t_2x$ are edges.
For each $C$-major vertex $v$ of $G\setminus N_3$, let $Q_1(v)$ be the $(v_1,v)$-gap in $Q$, with ends $q_1(v)$ and $q_1$, and let
$Q_2(v)$ be the $(x,v)$-gap in $Q$, with ends $q_2(v)$ and $q_2$. Note that $Q_1(v)$ is long, since $(v_1,B)$ is a remote base.
\\
\\
(1) {\em For every $C$-major vertex $v$ of $G\setminus N_3$,  either $v$ has a neighbour in $V(T)$, 
or $Q_2(v)$ has length less than $\ell$.}
\\
\\
Let $v$ be a $C$-major vertex of $G\setminus N_3$, and suppose that $v$ has no neighbour in $V(T)$, and 
$Q_2(v)$ is long. Suppose first that $v,v_1$ are adjacent. Now $T$
is odd, since $v_1$ has odd $x$-parity; and also $Q_2(v)$ is odd, since $v$ has odd $x$-parity. But $v_1$
has no neighbour in $Q_2(v)$, since $q_1$ is the only neighbour of $v_1$ in $Q$, and $q_1\notin V(Q_2(v))$ since $v,q_1$
are not adjacent. Since $T$ is anticomplete to $Q$ and hence to $Q_2(v)$,
the hole $v_1\dd v\dd q_2(v)\dd Q_2(v)\dd q_2\dd x\dd T^+\dd v_1$
is a long odd hole shorter than $C$, a contradiction. 

Now suppose that $v,v_1$ are nonadjacent.
Since $Q$, $Q_1(v)$ and $Q_2(v)$
are odd, and $Q_1(v)$, $Q_2(v)$ are edge-disjoint, it follows that the path $q_1(v)\dd Q\dd q_2(v)$ is odd. There is an induced path
$R$ between $q_1,q_2$ with interior in $V(T)\cup \{x,v_1\}$, and there are no edges between the interior of $R$
and the interior of $Q$, or between the interior of $R$ and $v$. Consequently if $q_1(v),q_2(v)$ are nonadjacent, the paths
$q_1(v)\dd v\dd q_2(v), q_1(v)\dd Q\dd q_2(v)$ and $q_1(v)\dd Q_1(v)\dd q_1\dd R\dd q_2\dd Q_2(v)\dd q_2(v)$ form a long jewel, and so 
there is a long odd hole of length less than $C$, a contradiction. Hence $q_1(v), q_2(v)$ are adjacent. Since the paths 
$Q_1(v)$ and $T$ are disjoint and anticomplete, there is a long pyramid formed by 
$v_1\dd q_1\dd Q_1(v)\dd q_1(v)$, a path between $v_1,q_2(v)$ with interior in $V(T)\cup \{x\}\cup V(Q_2(v))$
and a path between $v_1,v$ with interior in the interior of $D$, a contradiction. This proves (1).

\bigskip

It follows from (1) that there
is a set $\mathcal{F}_4$ of paths of $C$ with cost at most $2\ell-1$ that catches every $C$-major vertex of $G\setminus N_3$.
Then $\mathcal{F}_1\cup \mathcal{F}_2\cup \mathcal{F}_3\cup \mathcal{F}_4$ catches every $C$-major vertex of $G$ that is 
nonadjacent to $x$. This set has cost at most 
$$(5\ell+15)+ (4\ell+2)+ (5\ell+15)+ (2\ell-1)= 16\ell+28.$$ 
This proves \ref{cleannonnbrs}.~\bbox

\section{The algorithm}

In this section we will use the results of the previous sections to give our main theorem. We need
\begin{thm}\label{inconst0}
There is an algorithm with the following specifications:
\begin{description}
\item [Input:] A candidate $G$.
\item [Output:] Either decides that $G$ has a long odd hole, or decides that $G$ does not have a shortest long odd hole $C$ and
a remote base $(x,D)$ such that all $C$-major vertices are equal or adjacent to $x$.
\item [Running time:] $O(|G|^{9})$.
\end{description}
\end{thm}
\Proof
Suppose that $C$ is a shortest long odd hole with a remote base $(x,D)$, and every other $C$-major vertex is adjacent to $x$.
Here is an algorithm: we guess $x$ and the ends $d_1,d_2$ of the path $D$, and its middle vertex $d_0$ ($D$ has even length).
Let $Z$ be the set consisting of $x$ and all its neighbours. 
For $i = 1,2$, compute the set $D_i$ of internal vertices of all shortest paths with interior in $G\setminus Z$ 
between $d_i,d_0$. 
Let $Y$ be the set of vertices of $G$
not in $D_1\cup D_2\cup \{d_0,d_1,d_2\}$ and with a neighbour in $D_1\cup D_2\cup \{d_0\}$. 
Apply \ref{testcleanhole} to $G\setminus (Y\cup \{x\})$. If it finds a long odd hole, output that 
$G$ has a long odd hole. If (after checking all possible guesses) we did not find a long odd hole, output that
$G$ has no shortest long odd hole $C$ and
base $(x,D)$ such that every other $C$-major vertex is adjacent to $x$.

To see correctness, we only need to check correctness when $G$ has a shortest long odd hole $C$ and a remote base $(x,D)$, and
every other $C$-major vertex is adjacent to $x$. In this case, when we guess correctly,  
the interior of $D$ contains no vertices in $Z$; and since the subpath of $D$
between $d_1,d_0$ has length less than $|C|/2$, it is a shortest path of $G$ between $d_1,d_0$ containing no $C$-major vertices,
by \ref{shortpath}.
Any shortest path with interior in $G\setminus Z$
between $d_1,d_0$ will contain no $C$-major vertices, because they all belong to $Z$. Hence the interior of $D$
is a subset of the set $D_1\cup D_2\cup \{d_0\}$ computed by the algorithm. Also by \ref{shortpath}, $V(C)\setminus V(D)$
is anticomplete to $D_1\cup D_2\cup \{d_0\}$; and so the set $Y$ computed by the algorithm contains no vertices of $C$.
But it does contain all the $C$-major vertices except $x$; because they all have neighbours in
the interior of $D$, and do not belong to $D_1\cup D_2$ since $D_1\cup D_2$ is disjoint from $Z$. Hence after deleting
$Y\cup \{x\}$, all $C$-major vertices have been deleted and $C$ has become a clean shortest long odd hole, and the algorithm
of  \ref{testcleanhole} will detect a long odd hole; and so the output is correct.

For running time, there are $|G|^{4}$ guesses to check, and each one takes time $O(|G|^5)$ (because we are applying
 \ref{testcleanhole}). Thus the total running time is $O(|G|^{9})$. This proves \ref{inconst0}.~\bbox

Now our main result \ref{mainthm}, which we restate:

\begin{thm}\label{mainthm2}
There is an algorithm with the following specifications:
\begin{description}
\item [Input:] A graph $G$.
\item [Output:] Decides whether $G$ has a long odd hole.
\item [Running time:] $O(|G|^{20\ell+40})$.
\end{description}
\end{thm}
\Proof
First we apply \ref{testcandidate}, and we may assume we determine that $G$ is a candidate.
Next we apply \ref{testcleanhole} to $G$, and we assume we find that there is no clean shortest long odd hole.
So either $G$ has no long odd hole, or 
it has a shortest long odd hole $C$ with a base $(x,D)$.
We assume the latter. Let $R_1,R_2$ be the paths of $C$ of length $2\ell$
with middle vertices the ends $d_1,d_2$ of $D$. We guess $x$ and $R_1,R_2$, and delete all vertices caught by $\{R_1,R_2\}$ 
different from $x$, producing $H$ say, and run \ref{testcleanhole} on $H\setminus \{x\}$; and assuming it still does not 
find a long odd hole, and $G$ has a long odd hole, then when we guess correctly, $H$ has a shortest long odd hole $C$
with a remote base $(x,D)$.

Now by \ref{cleannonnbrs} there is a set $\mathcal{F}$ of paths of $C$, with cost at most $16\ell+28$, that catches every $C$-major
vertex of $H$ nonadjacent to $x$. We guess the paths in $\mathcal{F}$, and delete all the vertices they catch. Then we apply \ref{inconst0}
to the resulting graph. If it finds a long odd hole we are done; and it will do so in the case when we guess correctly.
If after all guesses we never find a long odd hole, we return that there is none. 

The total cost of guesses for the paths is $(4\ell+2)+(16\ell+28)=20\ell+30$, and we also have to guess $x$; and checking
each guess takes time $O(|G|^9)$, since we are applying \ref{inconst0}. Thus the total running time is $O(|G|^{20\ell+40})$.
This proves \ref{mainthm2}.~\bbox

\section*{Acknowledgement}
We would like to thank Hou Teng Cheong, who found a significant error in an earlier version of this paper.

\end{document}